\documentclass[a4paper,10pt]{article}

\usepackage{amsmath,amsfonts,amssymb,latexsym,amsthm}
\title{On Pr\'ekopa-Leindler inequalities on metric-measure spaces}
\author{Erwan~\textsc{Hillion}}

  \theoremstyle{plain}
  \newtheorem{thm}{Theorem}[section]
  \newtheorem{cor}[thm]{Corollary}
  \newtheorem{prop}[thm]{Proposition}
  \newtheorem{lem}[thm]{Lemma}
  \newtheorem{defi}[thm]{Definition}

\DeclareMathOperator{\tr}{Tr}
\DeclareMathOperator{\hess}{Hess}
\DeclareMathOperator{\ric}{Ric}
\DeclareMathOperator{\diam}{Diam}

\begin{document}

\maketitle

\begin{abstract}
This work is devoted to the geometric analysis of metric-measure spaces satisfying a Pr\'ekopa-Leindler or a more general Borell-Brascamp-Lieb inequality.

Completing the early investigations by Cordero-Erausquin, McCann and Schmuckenschl\"{a}ger, we show that these functional inequalities characterize lower bounds on the Ricci curvature on a Riemannian manifold, providing thus an alternate version of Ricci curvature lower bounds in measured length spaces to the recent developments by Lott, Villani and Sturm.
We also investigate stability properties and geometric and functional inequalities, such as logarithmic Sobolev inequality and Bishop-Gromov diameter estimate, in measured length spaces satisfying a Pr\'ekopa-Leindler or a Borell-Brascamp-Lieb inequality.
\end{abstract}

\section{Introduction}

This paper is about analysis and geometry in measured length spaces, and is strongly related to the theory developed by Lott and Villani in \cite{LottVill}, and Sturm in \cite{Sturm06}, \cite{Sturm06II} in the recent years. In both these papers and the present work, the major goal is to extend the geometric notion of Ricci curvature lower bounds from the Riemannian setting to the more general setting of measured length spaces, and to generalize theorems allowing for geometric and analytic inequalities under those extended notions of Ricci curvature bounds.

A length space is a metric space $(X,d)$ such that the distance between two points $x,y \in X$ is the infimum of lengths of curves joining these points. A length space is called a geodesic space if, for every couple $x,y \in X$, this infimum is attained by at least one curve. These curves are called geodesics between $x$ and $y$. It is easy to see that any compact length space is a geodesic space. Basic theory of length spaces can be found in the monograph by Burago, Burago and Ivanov \cite{burago}.

A measured length space $(X,d,\nu)$ is a length space $(X,d)$ with a probability measure $\nu$ on the Borel $\sigma$-algebra of $(X,d)$. Any compact Riemannian manifold with its (normalized) Riemannian volume is a measured geodesic space. Another example is provided by considering $\mathbb{R}^n$ with the Euclidean distance and a Gaussian measure. More generally, weighted Riemannian manifolds provide a large class of measured length spaces; see \cite[Def.-prop. 2.91]{GHL}.

The class of (measured) metric spaces is usually endowed with the (measured) Gromov-Hausdorff topology: a sequence $(X_n,d_n,\nu_n)$ of measured compact length spaces converges to another measured length space $(X,d,\nu)$ in the sense of Gromov-Hausdorff if there exists a sequence of measurable maps $f_n : X_n \rightarrow X$, and a sequence $(\varepsilon_n)_n \rightarrow 0$ of positive numbers such that: \begin{equation*}
\begin{array}{rl}
 {\it (i)} & \forall x_n,y_n \in X_n, \ |d_n(x_n,y_n)-d(f_n(x_n),f_n(y_n))| < \varepsilon_n, \\
 {\it (ii)} & \forall x \in X \ \exists x_n \in X_n, \  d(f_n(x_n),x) < \varepsilon_n, \\
 {\it (iii)} & \displaystyle (f_n)_*\nu_n \underset{{n \rightarrow \infty}}{\rightharpoonup} \nu.
\end{array}
\end{equation*} See \cite[Chapter 7]{burago} or \cite[Chapter 27]{villani2008optimal} for a more general definition.

It is quite natural to ask if it is possible to extend geometric notions from the Riemannian setting to the setting of measured length spaces. For example, the notion of sectional curvature bounded below by a real number $K$ has been extended to length spaces by Alexandrov, by considering geodesic triangles. See \cite{burago} for more explanations.

The Ricci curvature is another important way to describe the curvature of a Riemannian manifold. As the sectional curvature, it can be deduced from the Riemann curvature tensor. The Ricci curvature tensor associates to each tangent space $T_m M$ of a Riemannian manifold $(M,g)$ a symmetric bilinear form $\text{Ric}_m$. A Riemannian manifold $(M,g)$ is said to have its Ricci curvature bounded below by a real number $K$ if, for every point $m \in M$, the inequality $\text{Ric}_m \geq K \textrm{Id}$ holds in the sense of symmetric bilinear forms on the vector space $T_m M$. A rigorous definition of Ricci curvature and a survey of its basics properties can be found in any textbook of Riemannian geometry, for example \cite[Def. 3.18]{GHL}.

The Ricci curvature describes the behavior of a volume element transported along a geodesic. Hence, the theorems that can be deduced from the assumption $\text{Ric}\geq K$ mainly deal with the behaviour of volume of balls (such as the Bishop-Gromov theorem, see \cite[Theorem 3.101]{GHL}), or, by the use of change of variables formula, with functional inequalities involving integrals, for example the logarithmic Sobolev inequality. For a short and understandable presentation of the Ricci curvature and a comprehensive list of the consequences of Ricci curvature bounds, the reader is referred to \cite[Chapter 14]{villani2008optimal}.

Ricci curvature type lower bounds have been extended from the Riemannian setting to the setting of measured length spaces by Lott and Villani in \cite{LottVill}, and independently but in a very similar way by Sturm in \cite{Sturm06}, \cite{Sturm06II}. In these papers, the authors first generalize the theory of optimal transportation to the setting of compact length spaces. This theory provides a distance $W_2$ on the space of $P(X)$ probability measures on a compact length space $(X,d)$; namely, $$W_2(\mu_0,\mu_1)^2 = \inf_{\pi \in \Pi(\mu_0,\mu_1)} \int_{X \times X} d(x_0,x_1)^2 d\pi(x_0,x_1),$$ where $\Pi(\mu_0,\mu_1)$ is the set of probability measures on $X\times X$ having $\mu_0$ and $\mu_1$ as first and second marginals. The metric space $(P(X),W_2)$ is also a length space, called the Wasserstein space. The study of the convexity of certain functionals along the geodesics of the Wasserstein space leads to the definition of a curvature-dimension condition $CD(K,N)$ on the metric space, where $K \in \mathbb{R}\cup\{-\infty\}$ and $N \in [1,\infty]$. 

For example, the $CD(K,\infty)$ condition holds if, for every Wasserstein geodesic $(\mu_t)_{t \in [0,1]}$ in the set $P_2^{ac}(X,d,\nu)$ of probability measures absolutely continuous w.r.t. the reference measure $\nu$, the entropy functional $t \mapsto H_\nu(\mu_t) := \int_X \rho_t log(\rho_t) d\nu$, defined for $t \in [0,1]$, is $K$-convex with respect to $W_2$ in the sense that $$H_\nu(\mu_t) \leq (1-t) H_\nu(\mu_0)+t H_\nu(\mu_1)-\frac{Kt(1-t)}{2} W_2(\mu_0,\mu_1)^2.$$ The definition of the condition $CD(K,N)$, for $N \neq \infty$ is similar, the entropy functional $H_\nu$ being replaced by the Renyi entropy $H_{N,\nu}(\mu) := N-N \int_X \rho^{1-\frac{1}{N}} d\nu.$ See \cite{LottVill} for the case where $K=0$, and \cite{Sturm06}, \cite{Sturm06II} for the general case.

In order to prove that the condition $CD(K,N)$ on a measured length space is a natural generalization of both conditions $\text{Ric} \geq K$ and $\textrm{dim}(M) \leq N$ on a Riemannian manifold, Lott, Villani and Sturm showed that it satisfies these three specifications:

\medskip

(i) \textbf{Coherence:} If $(X,d,e^{-V}d\nu)$ is a Riemannian manifold $M$ with its Riemannian distance and a weighted measure, it satisfies $CD(K,N)$ if and only if $\rm{\textrm{dim}}(M) \leq N$ and the Ricci curvature tensor satisfies $$ \forall x \in M, \ \text{Ric}_x+ \hess_x V \geq K g_x.$$

(ii) \textbf{Stability:} If $(X_n,d_n,\nu_n)$  is a sequence of compact measured length spaces converging, in the measured Gromov-Hausdorff sense, to a compact measured length space $(X,d,\nu)$, and if each $(X_n,d_n,\nu_n)$ satisfies $CD(K_n,N)$ for a converging sequence $K_n \rightarrow K$, then the limit space also satisfies $CD(K,N)$.

\medskip

(iii) \textbf{Usefulness:} The condition $CD(K,\infty)$, where $K>0$, yields a large class functional inequalities (e.g. Poincar\'e, log-Sobolev, transportation inequalities), that are already known to be consequences of the condition $\text{Ric}\geq K >0$ in the Riemannian setting. Moreover, the condition $CD(0,N)$, where $N \geq 1$ allows for a control of the growth of volume of balls, which is analogous to the Bishop-Gromov theorem in Riemannian geometry.

\medskip

An important step in the development by Lott, Villani and Sturm was accomplished in 2001 by Cordero-Erausquin, McCann and Schmuckenschlager in \cite{CMS01}: they namely proved that on a Riemannian manifold satisfying $\text{Ric}\geq0$, the entropy functional is convex along Wasserstein geodesics.
The main goal of \cite{CMS01}, and its successor \cite{CMS06}, was actually the study of a large class of functional inequalities on a Riemannian manifold, called the Borell-Brascamp-Lieb (or shortly $BBL$) inequalities under Ricci curvature bounds. As the $CD(K,N)$ condition, these inequalities depend on two parameters $K \in \mathbb{R}$ and $N \geq 1$. The present paper will focus only on two particular subclasses of the $BBL$ inequalities.

The first subclass of inequalities is the family of Pr\'ekopa-Leindler inequalities. It only depends on the parameter $K$. These inequalities can be seen as $BBL(K,N)$ inequalities where the parameter $N$ tends to the infinity. They are defined as follows.

\begin{defi}
A measured length space $(X,d,\nu)$ is said to satisfy the Pr\'ekopa-Leindler inequality with constant $K \in \mathbb{R}$, or shortly $PL(K)$, if, for every $t \in [0,1]$ and every triple of functions $(u,v,w) \in L^1(X,\mathbb{R}_+^*)$ satisfying $$\forall z \in Z_t(x,y), \ u(x)^{1-t} v(y)^t \exp\left( -K \frac{t(1-t)}{2}d(x,y)^2\right) \leq w(z),$$ the following inequality holds: $$\left(\int_X u d\nu\right)^{1-t} \cdot \left(\int_X v d\nu\right)^t \leq \int_X w d\nu.$$
\end{defi}

The set $Z_t(x,y)$ is the set of points $z \in X$ such that $d(x,z)=td(x,y)$ and $d(z,y)=(1-t)d(x,y)$. It is often called the set of $t$-intermediate points between $x$ and $y$. More generally, the set $Z_t(A,B)$ of $t$-intermediate points between two sets $A,B \subset X$ is defined by $Z_t(A,B) := \{Z_t(x,y) \ | \ x\in A, y \in B\}$.

The family of Pr\'ekopa-Leindler inequalities, especially $PL(0)$, have first been introduced in the Euclidean setting by Pr\'ekopa in \cite{Preko71}, and Leindler in \cite{Leindler} to generalize the multiplicative Brunn-Minkowski inequality, which has a more geometric content.

\begin{defi}
A measured length space $(X,d,\nu)$ is said to satisfy the mutliplicative Brunn-Minkowski inequality, or $BM$ if, for every measurable sets $A$ and $B$, and every $t$ in $[0,1]$, $$\nu^*(Z_t(A,B)) \geq \nu(A)^{1-t} \nu(B)^t,$$ where $\nu^*$ is the outer measure associated to $\nu$.
\end{defi}

It is easy to see that the $PL(0)$ condition implies the Brunn-Minkowski inequality, by taking indicator functions. However, Pr\'ekopa-Leindler inequalities are easier to manipulate that the Brunn-Minkowski inequality, and it is now simpler, in $\mathbb{R}^n$, to prove the $PL(0)$ inequality and then derive the Brunn-Minkowski inequality, than proving it directly.

\medskip

The second subclass of $BBL$ inequalities is obtained by taking the parameter $K$ equal to $0$. This is the family of $BBL(0,N)$ inequalities. Their definition requires to consider $p$-means between real numbers.

\begin{defi}
For $p \in \mathbb{R} \setminus \{0\}$, the $p$-mean between $a,b \geq 0$ is defined by: $\mathcal{M}_t^p(a,b) := ((1-t)a^p+tb^p)^{1/p}$ if $ab\neq 0$, and $\mathcal{M}_t^p(a,b) := 0$ if $ab=0$.\\
For $p=0$, define $\mathcal{M}_t^p(a,b) := a^{1-t}b^t$ if $ab\neq 0$, and $\mathcal{M}_t^p(a,b) := 0$ if $ab=0$.
\end{defi}

This $p$-mean is used to define a ``$p-BBL(0,N)$'' inequality.

\begin{defi}
Let $N \geq 1$ and $p \geq -1/N$ be real numbers. A measured length space $(X,d,\nu)$ is said to satisfy the $BBL_p(0,N)$ inequality if, for every $t\in [0,1]$ and every triple of functions $f,g,h : X \rightarrow \mathbb{R}_+$ satisfying $$\forall z \in Z_t(x,y), \ h(z) \geq \mathcal{M}_t^p(f(x),g(y)),$$ one has $$\int_X h d\nu \geq \mathcal{M}_t^{p/(1+Np)}\left(\int_X f d\nu, \int_X g d\nu \right).$$
\end{defi}

The Pr\'ekopa-Leindler inequality $PL(0)$ corresponds to the case where $p=0$ and $N=\infty$. The $BBL(0,N)$ condition is the intersection of a family of $BBL_p(0,N)$ conditions.

\begin{defi}
A measured length space $(X,d,\nu)$ is said to satisfy the $BBL(0,N)$ inequality for a certain $N \geq 1$ if it satisfies $BBL_p(0,N)$ for every $p \geq -1/N$.
\end{defi}

In \cite{CMS01}, the authors proved that the $BBL_p(0,N)$ inequality implies the $BBL_{p'}(0,N)$ inequality, for any $p'\geq p$. This observation yields in particular that any Riemannian manifold satisfying the $BBL(0,N)$ condition for a given $N\geq 1$ also satisfies the $PL(0)$ condition. This fact will be useful in the proof of Theorem \ref{BBLRic}.

\medskip

In the same paper, the authors also proved  that the $PL(K)$ condition is satisfied in any Riemannian manifold whose Ricci curvature tensor satisfies $\text{Ric} \geq K$, and that the $BBL(0,N)$ condition is satisfied on any Riemannian manifold with $\text{Ric} \geq 0$ and $\textrm{dim}(M) \leq N$.
Actually, the authors went a bit further by defining $BBL(K,N)$ inequalities on a Riemannian manifold for a general $K \in \mathbb{R}$. The extension of theorems proven in this paper under a $PL(K)$ or a $BBL(0,N)$ condition to the more general $BBL(K,N)$ condition still remains an open question.

\medskip

General $BBL(K,N)$ inequalities have been recently considered by Bacher in \cite{Bacher} in a length space setting. In this paper, the author proved the implication $CD(K,N) \Rightarrow BBL(K,N)$ in the case of compact non-branching metric spaces, and proved the stability of these inequalities by Gromov-Hausdorff convergence of metric spaces.

\medskip

So far, it seems that the natural specifications (i),(ii) and (iii), satisfied by the $CD(K,\infty)$ and $CD(0,N)$ conditions, are also partially satisfied by the conditions $PL(K)$ and $BBL(0,N)$: the whole point (ii), about stability of the inequalities by Gromov-Hausdorff convergence, is the main result of Bacher's article (\cite[Theorem 4.11]{Bacher}). The ``if'' part of the first point, about coherence with the Riemannian setting, has been proven in \cite{CMS01}, and in \cite{CMS06} for the weighted Riemannian case. Finally, the $PL(K)$ condition, with $K>0$, implies functional inequalities, at least in a smooth setting; see for instance \cite{BobkovLedoux} for a proof of the logarithmic Sobolev inequality in the Euclidean case.

\medskip

The main goal of this article is to complete the picture that the natural specifications (i), (ii) and (iii), satisfied by the $CD(K,\infty)$ condition, are also fully satisfied by the condition $PL(K)$, and that the condition $BBL(0,N)$ is a natural generalization of both conditions $\text{Ric} \geq 0$ and $\textrm{dim} \leq N$ on a Riemannian manifold.

\medskip

The first part of this paper focus on the study of the $PL(K)$ condition on a measured length space. We first prove, in Theorem \ref{PLRic}, that on a weighted Riemannian manifold $(M,g,e^{-V}d\textrm{vol})$ on which the $PL(K)$ condition is satisfied, one has necessarily $\text{Ric} + \hess V\geq K$, which corresponds to the ``only if'' part of point (i). Then we give, in Theorem \ref{PLTenso}, a result about the tensorization properties of the condition $PL(K)$: for a suitable product distance, if the condition $PL(K)$ condition is satisfied on two measured length spaces $(X_1,d_1,\nu_1)$ and $(X_2,d_2,\nu_2)$, then it is also satisfied on the product $(X_1 \times X_2, d_{X_1 \times X_2}, \nu_1 \otimes \nu2)$. We finally give some analytic consequences of the $PL(K)$ condition, for $K>0$, namely the logarithmic Sobolev (Theorem \ref{LogSob}) and Poincar\'e (Proposition \ref{Poinc}) inequalities, and Talagrand's transportation inequality (Theorem \ref{Talagrand}).

\medskip

The second part of the paper is about the study of the $BBL(0,N)$ condition on a measured length space, and follows the same lines as the first part. We first focus on the Riemannian case, proving the equivalence between the $BBL(0,N)$ condition and the couples of conditions $\text{Ric} \geq 0$ and $\textrm{dim} \leq N$ on a Riemannian manifold (Theorem \ref{BBLRic}).We then prove a tensorization property (Theorem \ref{BBLTenso}), which enlightens the intuitive meaning of the parameter $N$ as an upper bound on the dimension of the underlying metric space. Finally, we show that the $BBL(0,N)$ condition allows for the control of the growth of volume of balls centered around a fixed point (Theorem \ref{BishopGromov}), which generalizes in the length space setting a particular case of the Bishop-Gromov theorem in Riemannian geometry. As a consequence of Theorem \ref{BishopGromov}, it becomes possible to control the Hausdorff dimension of a measured length space satisfying the $BBL(0,N)$ condition by $N$ (Corollary \ref{HausDim}). Another consequence of Theorem \ref{BishopGromov} and Theorem \ref{Talagrand} is Theorem \ref{BonnetMyers}, that allows for the control of the diameter of a measured length space satisfying both conditions $BBL(0,N)$ and $PL(K)$, for $N \geq 1$ and $K>0$, and which is a weak version of the Bonnet-Myers Theorem in Riemannian geometry.

\section{The Pr\'ekopa-Leindler inequality}

This section is devoted to the study of Pr\'ekopa-Leindler inequalities on a measured length space. We first link the $PL(K)$ condition on a Riemannian manifold to uniform lower bounds of the Ricci curvature tensor (Theorem \ref{PLRic}). We then prove a tensorization theorem that allows us to deduce the $PL(K)$ inequality on a tensor product of two measured length spaces, with a suitable product distance, from the $PL(K)$ inequality on these spaces (Theorem \ref{PLTenso}). We finally show that the $PL(K)$ inequality, with $K>0$, is a sufficient condition for a (compact) measured length space to have functional inequalities, such as the logarithmic Sobolev inequality (Theorem \ref{LogSob}), the Poincar\'e inequality (Theorem \ref{Poinc}), and the transportation inequality $T_2$ (Theorem \ref{Talagrand}).

\subsection{The Riemannian case}
In this paragraph, we prove the converse of the main theorem stated and proven by Cordero-Erausquin, McCann and Schmuckenschl\"{a}ger in \cite{CMS01}. This theorem asserts that, on a weighted Riemmannian manifold $(M,g,\nu)$, with $\nu := e^{-V}d \textrm{vol}$ and $\hess V + \text{Ric} \geq K$, the Pr\'ekopa-Leindler inequality $PL(K)$ holds.
We establish here the converse statement: 

\begin{thm}\label{PLRic}
Let $(M,g,e^{-V}d\textrm{vol})$ be a weighted Riemannian manifold satisfying the Pr\'{e}kopa-Leindler inequality with constant $K$. Then the Ricci curvature tensor satisfies $$\ric_m + \hess_m V \geq K$$ at each point $m$ of the manifold.
\end{thm}

As a consequence, this theorem and its converse imply that, on a weighted Riemannian manifold $(M,g,e^{-V}d\textrm{vol})$, both conditions $\ric + \hess V \geq K$ and $PL(K)$ are equivalent. Moreover, the stability results \cite[Theorem 5.19]{LottVill} and \cite[Theorem 4.11]{Bacher} imply that the conditions $PL(K)$ and $CD(K,\infty)$ are equivalent on Gromov-Hausdorff limits of weighted Riemannian manifolds. In the general setting of measured length spaces, the equivalence between the conditions $PL(K)$ and $CD(K,\infty)$ is still an open question. See (\cite[Theorem 3.1]{Bacher}) for the implication $CD(K,\infty) \Rightarrow PL(K)$ in the case where the underlying length space is non-branching.\\

The proof of Theorem \ref{PLRic} follows the same lines as the proof of its converse, and the reader is referred to \cite{CMS06} for a better understanding of the methods. It is also inspired by the proof of the implication $CD(K,\infty) \Rightarrow \ric \geq K$ given by Lott and Villani in \cite[Theorem 7.3]{LottVill}. The starting point is the following existence result.

\begin{prop}\label{theta-existence}
Let $(M,g)$ be a Riemannian manifold. Given a point $m \in M$ and a tangent vector $\tilde{v} \in T_m M$, there exists a smooth real function $\theta$ defined on an open ball $B(m,r_0)$ such that:

(i) $v := \nabla_m \theta$ is colinear to $\tilde{v}$.

(ii) $\hess_m \theta = 0$

(iii) $\forall x \in B(m,r_0), \ \gamma_x : t \mapsto \exp_x(t \nabla_x \theta)$ defines a minimal geodesic.

(iv) $\forall t \in [0,1], \ F_t : x \mapsto \exp_x(t \nabla_x \theta)$ is a diffeomorphism from $B(m,r_0)$ onto its image.
\end{prop}

\textbf{Proof:} The existence of a function $\tilde{\theta}$ satisfying both points (i) and (ii) is a classical result which can be obtained by using local coordinates. See \cite[Theorem 14.8]{villani2008optimal} for another use of this result. Point(iii) comes from a classical result of Riemannian geometry: changing the function $\tilde{\theta}$ into a proportional function $\theta := \eta \tilde{\theta}$, it is possible to apply \cite[Theorem 2.92]{GHL}, which gives the result. The last point is easy to prove in local coordinates, by applying the inverse function Theorem. $\square$

\medskip

We will also use a classical result about invertible matrices of Jacobi fields.

\begin{lem}\label{jacobi-preko}
Let $A(t)$ be a matrix of invertible Jacobi fields along a geodesic $\gamma$, written in an orthonormal moving frame along $\gamma$. 
If $\phi : [0,1] \rightarrow \mathbb{R}$ is defined by $\phi(t) := -\log(\det(A(t)))$, then: $$\phi''(t) = \ric_{\gamma(t)}(\gamma'(t)) + \tr((A'(t)A(t)^{-1})^2).$$
\end{lem}

\textbf{Proof:}
Recall that the function $t \mapsto A(t)$ satisfies the Jacobi equation: $$A''(t) = R(t) A(t),$$ where $R(t)$ is a matrix whose trace is equal to $\ric_{\gamma(t)}(\gamma'(t))$.
Set $B(t) := A'(t)A(t)^{-1}$. Simple calculations yield
\begin{eqnarray*}
\phi'(t) &=& -det(A(t)) \tr(A'(t)A(t)^{-1}) det(A(t))^{-1} \\
&=& \tr(B(t)),
\end{eqnarray*} and \begin{eqnarray*}
\phi''(t) &=& \tr(B'(t)) \\ &=& \tr(A''(t)A(t)^{-1})+\tr(A'(t)A'(t)A(t)^{-2}) \\ &=& \tr(R(t))+\tr(B(t)^2), 
\end{eqnarray*} and this is the expected result. $\square$

\bigskip

\textbf{Proof of Theorem \ref{PLRic}:} Let $m \in M$, and $\tilde{v} \in T_m M$. 
We want to prove that $\text{Ric}_m(v,v)+(\hess_m V)(v,v) \geq K |v|^2$ for a certain vector $v \in T_m M$ colinear to $\tilde{v}$.
Let $\theta : B(m,r_0) \rightarrow M$ be a function satisfying the conclusion of proposition ~\ref{theta-existence}. We denote by $v$ the vector $\nabla_m \theta$. The vector $v$ is colinear to $\tilde{v}$ by (i) of Proposition~\ref{theta-existence}.
Let $r_0>r>0$ and $\mu_0$ be a probability measure having a smooth density $\rho_0$ supported in $B(m,r)$.
Set $$F_t(x) := \exp(t \nabla_x \theta) =: \gamma_x(t),$$ and $$\mu_t := (F_t)_*\mu_0.$$ 

By (iii) of Proposition \ref{theta-existence}, each $\gamma_x$ is a minimal geodesic of the manifold. 
Each $\mu_t$ is a probability measure on $M$. Denote by $\rho_t$ the density of $\mu_t$ with respect to the Riemannian volume element $d\textrm{vol}$. The densities $\rho_t$ are ruled by the Monge-Amp\`{e}re equation, \begin{equation}\label{PLRic-eq1}\rho_0(x) = \rho_t(\gamma_x(t)) \det(d(F_t)(x)).\end{equation} 
It is also known that, when $x$ is fixed, $t \mapsto A_x(t)$ is a matrix of Jacobi fields along the geodesic $\gamma_x$, where $A_x(t)$ is the matrix of $d(F_t)(x)$ in an orthonormal moving frame along $\gamma_x$. By (iv) of Proposition \ref{theta-existence}, each matrix $A_x(t)$ is invertible. Moreover, this is the only matrix of Jacobi fields along $\gamma_x$ satisfying the initial conditions $$A_x(0) = Id, \ A'_x(0)=0.$$
Writing $\phi_x(t) := -\log(\det(A_x(t)))$, Lemma \ref{jacobi-preko} gives \begin{equation}\label{PLRic-eq2}\phi_x''(t) = \ric_{\gamma_x(t)}(\gamma_x'(t)) + \tr((A_x'(t)A_x(t)^{-1})^2).\end{equation}
Taking $x=m$ and $t=0$, (\ref{PLRic-eq2}) becomes \begin{equation}\label{PLRic-eq3} \ric_m(v)=\phi_m''(0).\end{equation}
With these notations, another formulation of (\ref{PLRic-eq1}) is thus possible, namely \begin{equation}\label{PLRic-eq4}\rho_x(\gamma_x(t)) = \rho_0(x) \exp(\phi_x(t)).\end{equation}
Now, use the Pr\'ekopa-Leindler inequality, being careful of dealing with densities with respect to the reference measure $d\nu = \exp(-V)d\textrm{vol}$. More precisely, if $\mu_t=\rho_t d\textrm{vol} = u_t d\nu$, then $u_t(x) = \exp(V(x)) \rho_t(x)$. Let $t$ be in $[0,1]$ and $ h : F_t(B(m,r)) \rightarrow \mathbb{R}_+$ be defined by $$h(\gamma_x(t)) \exp(V(\gamma_x(t))) :=$$$$ u_0(\gamma_x(0))^{1-t} u_1(\gamma_x(1))^t \exp\left( - K \frac{t(1-t)}{2} d(\gamma_x(0),\gamma_x(1))^2\right).$$ This definition is consistent since, when $t$ is fixed, $F_t$ is a diffeomorphism from $B(m,r)$ onto $F_t(B(m,r))$. The curve $\gamma_x$ being a minimal geodesic implies: $d(\gamma_x(0),\gamma_x(1))=|\nabla \theta(x)|$.
As $u_0$ and $uo_1$ are probability densities w.r.t. $\nu$, the Pr\'ekopa-Leindler inequality yields $$\int_M h d\textrm{vol} = \int_M h \exp(V) d\nu \geq 1 = \int_M \rho_s(\gamma_s(x)) d\textrm{vol}.$$
On the other hand, (\ref{PLRic-eq4}) expresses that $$
\int_M \rho_0(x) \exp\left( (1-t) V(\gamma_x(0))+tV(\gamma_x(1)) \right) $$$$
\exp\left(-V(\gamma_x(t))+t \phi_x(1)-K \frac{t (1-t)}{2} |\nabla \theta(x)|^2 \right) d\nu $$$$
 \geq \int_M \rho_0(x) \exp(\phi_x(t)) d\nu.
$$
The above formula does not depend on the size of the ball $B(m,r)$ which contains the support of $\rho_0$. In particular, letting $r \rightarrow 0$, $$V(\gamma_m(t))+\phi_m(t) \leq (1-t)V(\gamma_m(0))+tV(\gamma_m(1))+t\phi_m(1)-K \frac{t (1-t)}{2} v^2.$$
This means that the function $t \mapsto V(\gamma_m(t))+\phi_m(t)$ is $K v^2$-convex. In particular, its second derivative at $t=0$ is greater than $K v^2$, and by (\ref{PLRic-eq3}) this yields the formula.$\square$

\subsection{Tensorization}

The simplest way to prove that the Euclidean space $\mathbb{R}^n$ satisfies the $PL(0)$ condition is to prove it first for $n=1$, and to use induction on $n$ to complete the general case. See \cite{GardnerBM} for a proof using this strategy. What makes the induction step work can be seen as a tensorization property: the condition $PL(0)$ on both $\mathbb{R}^{n_1}$ and $\mathbb{R}^{n_2}$, with Euclidean distance and Lebesgue measure, implies $PL(0)$ on their product space $\mathbb{R}^{n_1+n_2}$, with Euclidean distance and Lebesgue measure, which is the tensor product of Lebesgue measures on $\mathbb{R}^{n_1}$ and $\mathbb{R}^{n_2}$.\\The proof of this tensorization property let us hope it can be generalized for a general $PL(K)$ condition, and in others settings that the Euclidean one. Another clue for a possible generalization is the fact that some functional inequalities implied by the condition $PL(K)$ for $K>0$ (in particular the logarithmic Sobolev inequality, see Theorem \ref{LogSob}) are known to be stable by tensorization (ref : Gross75).
In this paragraph, we show that the Pr\'ekopa-Leindler inequality tensorizes for a certain distance on the product space, which can be seen as an ``Euclidean product distance'':

\begin{defi}
The Euclidean product distance on $X_1 \times X_2$ is defined by $$\forall \ ((x_1,x_2),(y_1,y_2)) \in (X_1 \times X_2)^2$$$$ \ d_{X_1\times X_2}((x_1,x_2),(y_1,y_2))^2 := d_{X_1}(x_1,y_1)^2+d_{X_2}(x_2,y_2)^2.$$
\end{defi}

\begin{thm}\label{PLTenso}
Let $(X_1,d_1,\nu_1)$ and $(X_2,d_2,\nu_2)$ be two measured length spaces and $K \in \mathbb{R}$. We suppose that both spaces satisfy $PL(K)$. Then the product space $(X_1\times X_2, d_{X_1 \times X_2}, \nu_1 \otimes \nu_2)$ also satisfies $PL(K)$.
\end{thm}

\textbf{Proof:} Denote by $(X,d)$ the metric space $(X_1\times X_2, d_{X_1 \times X_2})$.
Let $t$  be in $[0,1]$ and $u,v,w : X \rightarrow \mathbb{R}_+$ be three functions satisfying $$\forall x,y \in X \ , \ \forall z \in Z_t(x,y) \ , \ u(x)^{1-t} v(y)^t \exp\left( - K \frac{t(1-t)}{2} d(x,y)^2 \right) \leq w(z).$$
Fix $x_2,y_2,z_2 \in X_2$ and consider the functions $u_{x_2} := u(.,x_2), v_{y_2} := v(.,y_2), w_{z_2} := w(.,z_2) : X_1 \rightarrow \mathbb{R}_+$ . Let $x_1, y_1 \in X_1$ and $z_1 \in Z_t(x_1,y_1)$. Then $(z_1,z_2) \in Z_t((x_1,x_2),(y_1,y_2))$, so that $$\forall x_1,y_1 \in X_1 \ , \ \forall z_1 \in Z_t(x_1,y_1) ,$$ \begin{equation} u_{x_2}(x_1)^{1-t} v_{y_2}(y_1)^t \exp\left(-K\frac{t(1-t)}{2}(d_1(x_1,y_1)^2+d_2(x_2,y_2)^2)\right) \leq w_{z_2}(z_1).\end{equation} 
The triple of functions $$u_{x_2}, \ v_{y_2}, \ \exp\left(K \frac{t(1-t)}{2} d_2(x_2,y_2)^2 \right) w_{z_2}$$ satisfy the hypothesis for $PL(K)$ on the space $(X_1,d_1,\nu_1)$. Therefore,  $$\left( \int_{X_1} u(.,x_2) d\nu_1 \right)^{1-t} \cdot \left( \int_{X_1} v(.,y_2) d\nu_1 \right)^t $$ \begin{equation}\leq \exp\left(K \frac{t(1-t)}{2} d_2(x_2,y_2)^2 \right) \int_{X_1} w(.,z_2) d\nu_1.\end{equation}
Defining the function $\tilde{u} : X_2 \rightarrow \mathbb{R}$ by $$\tilde{u}(x_2) := \int_{X_1} u(.,x_2) d\nu_1,$$ and the functions $\tilde{v}, \tilde{w}$ in the same way, yields: $$\forall x_2,y_2 \in X_2, \ \forall z_2 \in Z_t(x_2,y_2),$$ \begin{equation}\tilde{u}(x_2)^{1-t} \tilde{v}(y_2)^t \exp\left(-K \frac{t(1-t)}{2} d_2(x_2,y_2)^2\right) \leq \tilde{w}(z_2).\end{equation} The triple of functions $(\tilde{u},\tilde{v},\tilde{w})$ satisfy the hypothesis of the $PL(K)$ inequality on $(X_2,d_2,\nu_2)$, so that $$\left(\int_{X_2} \tilde{u} d\nu_2\right)^{1-t} \left(\int_{X_2} \tilde{v} d\nu_2 \right)^t \leq \int_{X_2} \tilde{w} d\nu_2.$$
By immediate properties of tensor products, $$\left( \int_X u d\nu\right)^{1-t} \left( \int_X v d\nu \right)^t \leq \int_X w d\nu,$$ and this is the announced result. $\square$

\medskip

This proof can easily be adapted to the more general case where $(X_1,d_1,\nu_1)$ (resp. $(X_2,d_2,\nu_2)$) satisfy $PL(K_1)$ (resp. $PL(K_2)$). Then the product space satisfy $PL(K)$, where $K:=min(K_1,K_2)$.

Also remark that, by Theorem \ref{PLRic} and its converse, Theorem \ref{PLTenso} yields that the condition $\text{Ric} \geq K$ is preserved by product of Riemannian manifolds (see \cite[Definition 2.15]{GHL} for the definition).

\subsection{Analytic consequences of the $PL(K)$ condition}

We end this section by presenting some consequences that arise from the hypothesis that a measured length space satisfy the $PL(K)$ condition, for $K>0$. Most of these consequences are already known in smoother cases (i.e. in an Euclidean or Riemannian setting), and proofs are very similar. For more details, the reader is referred to references given before each proposition.

It is also worth observing that the consequences of the condition $PL(K)$ are also very similar to the consequences of the condition $CD(K,\infty)$; in particular, the constants found in functional inequalities are the same.

\subsubsection{Logarithmic Sobolev and Poincar\'e inequalities}

The logarithmic Sobolev inequality allows to control the entropy of a function by its $L^2$ norm. It appears as an infinite-dimensional and degenerate version of the Sobolev inequalities, which control a $L^p$ norm of a function by the $L^2$ norms of the function and of its gradient.
Such a control does not always exists out of the Riemannian setting. For example, on $\mathbb{R}^n$ with the standard Gaussian measure $\gamma_n$, the Sobolev inequality degenerates into $$H_{\gamma_n}(f^2) \leq 2 ||\nabla f||_2^2,$$ for $f$ enough regular, where $H_{\gamma_n}(f^2)$ is the relative entropy of the function $f^2$: $$H_{\gamma_n}(f^2) := \int_{\mathbb{R}^n}f^2 \log(f^2) \gamma_n(dx)-\int_{\mathbb{R}^n}f^2\gamma_n(dx) \cdot log \int_{\mathbb{R}^n}f^2\gamma_n(dx).$$ A smooth function $f$ whose gradient is in $L^2(\gamma_n)$ do not necessarily lie in a $L^p(\gamma_n)$ space with $p>2$, but only in the space $L^2logL^2$. This kind of inequality is called logarithmic Sobolev inequality.

\medskip

These inequalities have been widely studied since the first works of Nelson and Gross in the 1960-70's, and there is now a huge litterature on this subject. The reader is referred to \cite{ABCFGMRS} (in French) or \cite{Bakry97}.
Logarithmic Sobolev inequalities are used to find exponential bounds on the rate of convergence of Markov diffusions to their equilibrium. See \cite{BakrySF} for more details, and for more properties strongly linked with the logarithmic Sobolev inequality.\\
In \cite{BakryEmery84} and \cite{BakryEmery85}, Bakry and Emery, studying Markov diffusions, noticed the influence of Ricci curvature bounds on the constant of the logarithmic Sobolev inequality. In particular, their celebrated Bakry-Emery criterion yields the following.

\begin{prop}\label{BakryEmery}
Let $(M,g)$ be a Riemannian manifold with a reference probability measure $d\nu(x) = exp(-V(x)) d\textrm{vol}$. Suppose that there exist a number $K>0$ such that $$\ric_m+\hess_m V \geq K$$ at each point $m$ of the manifold. Then, for every function $f$ enough regular, $$H_\nu(f^2) \leq \frac{2}{K} \int_X |\nabla f|^2(x) d\nu(x).$$
\end{prop}

In this paragraph, we show that the conclusion of \ref{BakryEmery} holds under the assumption that the condition $PL(K)$ is satisfied. By Theorem \ref{PLRic}, it is immediate to notice this in the weighted Riemannian case. To prove it in the general measured length space setting, we first need to generalize the notion of gradient; in fact, this is only the notion of gradient norm wich is generalized:

\begin{defi}
Let $g: X \rightarrow \mathbb{R}$ be a continuous function on a length space $(X,d)$. The gradient norm of $g$ at a point $x \in X$ is defined by $$|\nabla g|(x) = \limsup_{z \rightarrow x} \frac{|g(x)-g(z)|}{d(x,z)}.$$
\end{defi}

There is a chain rule associated to the gradient norm:

\begin{prop}
Let $g: X \rightarrow \mathbb{R}$ be continuous and $h: \mathbb{R} \rightarrow \mathbb{R}$ be differentiable. Then $$|\nabla (h \circ g)|(x) = h'(g(x)) |\nabla g|(x).$$
\end{prop}

\medskip

We can now state the main result of this paragraph :

\begin{thm}\label{LogSob}
Let $(X,d,\nu)$ be a compact measured length space satisfying the Pr\'ekopa-Leindler inequality with a constant $K>0$. Then for every continuous function $f: X \rightarrow \mathbb{R}_+$: \begin{equation}\label{LogSobEq}H_{\nu}(f^2) \leq \frac{2}{K} \int_X |\nabla f|^2(x) d\nu(x).\end{equation}
\end{thm}

Note that the same conclusion, with the same constant $2/K$ can be deduced from the condition $CD(K,\infty)$, with $K>0$. This has been done in \cite[Theorem 6.1]{LottVill}.

\bigskip

\textbf{Proof:} Define the function $g: X \rightarrow \mathbb{R}$ by $f^2=\exp(g)$. Since $$|\nabla f|^2 = \left| \nabla \exp\left(\frac{g}{2}\right) \right| ^2 = \left(\frac{1}{2}|\nabla g| \exp\left(\frac{g}{2}\right)\right)^2 = \frac{1}{4} |\nabla  g|^2 \exp(g),$$ it is sufficient to prove that \begin{equation}H_{\nu}(\exp(g)) \leq \frac{1}{2K} \int_X |\nabla g|^2(x) \exp(g(x)) d\nu(x).\end{equation}
Writing, for $r>0$, $$|\nabla_r g|(x):= \sup_{z \in B(x,r)} \frac{|g(x)-g(z)|}{d(x,z)},$$ it is sufficient, by a Beppo-Levi argument, to prove that, for every $r>0$ small enough, \begin{equation}H_{\nu}(\exp(g)) \leq \frac{1}{2K} \int_X |\nabla_r g|^2(x) \exp(g(x)) d\nu(x).\end{equation}
Recall now the following formula that expresses the entropy of a function as a derivative of its $L^p$-norms: \begin{equation}H_\nu(f) = - \lim_{t \rightarrow 1} \frac{\int_X f d\nu - \left(\int_X f^{1/t}\right)^t}{1-t}.\end{equation}
Fix $t \in ]0,1]$, and apply the Pr\'ekopa-Leindler inequality with $u(x)=1$, $v(y)=\exp\left(\frac{g}{t}\right)$ and $w(z)=\exp(g_t(z))$, where $g_t$ is defined by $$g_t(z) = \sup_{z \in Z_t(x,y)}\left\{g(y)-K\frac{t(1-t)}{2}d(x,y)^2\right\}.$$
The $PL(K)$ inequality then expresses that \begin{equation}\left(\int_X \exp\left(\frac{g}{t}\right) d\nu \right)^t \leq \int_X \exp(g_t) d\nu.\end{equation}
Hence: \begin{equation}H_\nu(\exp(g)) \leq \limsup_{t \rightarrow 1} \int_X \frac{\left|\exp(g(z))-\exp(g_t(z))\right|}{1-t} d\nu(z).\end{equation}
Fixing $z \in X$, we have:

\begin{eqnarray}\frac{g_t(z)-g(z)}{1-t} &=& \sup_{(x,y):z\in Z_t(x,y)} \left\{\frac{g(y)-g(z)}{1-t}-K\frac{t}{2}d(x,y)^2 \right\}\\
&=& \sup_{(x,y):z\in Z_t(x,y)} \left\{\frac{g(y)-g(z)}{d(y,z)}d(x,y)-K\frac{t}{2}d(x,y)^2 \right\}. 
\end{eqnarray}

As $X$ is compact, its diameter $D$ is finite. Note that, if $z \in Z_t(x,y)$, then $d(y,z) = (1-t)d(x,y) \leq (1-t)D$, so $y \in B(z,(1-t)D)$, and: \begin{equation}\frac{g_t(z)-g(z)}{1-t} \leq \sup_{(x,y):z\in Z_t(x,y)} \left\{|\nabla_{(1-t)D} g|(z) d(x,y)-K\frac{t}{2}d(x,y)^2 \right\}.\end{equation}
Moreover, it is well known that the function $u \mapsto au-bu^2$, where $b>0$, attains its maximum at $u=\frac{a}{2b}$ and this maximum equals $\frac{a^2}{4b}$. Thus: \begin{equation}\frac{g_t(z)-g(z)}{1-t} \leq \frac{t}{2K} |\nabla_{(1-t)D} g|^2(z).\end{equation}
Letting $t$ go to $1$ completes the proof. $\square$ 

\medskip

It is well known that the logarithmic Sobolev inequality implies the Poincar\'e inequality, which controls the $L^2$ norm of a function by the $L^2$ norm of its gradient. The proof of this fact in a measured length space setting can be found in \cite[Theorem 6.18]{LottVill}. 

\begin{prop}\label{Poinc}
Let $(X,d,\nu)$ be a measured length space satisfying the logarithmic Sobolev inequality $LSI(K)$, with $K>0$. Then, for every Lipschitz function $h$, with $\int_X h d\nu =0$, \begin{equation}\int_X h^2 d\nu \leq \frac{1}{K} \int_X |\nabla h|^2 d\nu.\end{equation}
\end{prop}

Actually, a stronger fact has been proven by Lott and Villani $CD(K,N)$ in \cite{LottVillweak}: under a $CD(K,N)$ condition with $K>0$, the Poincar\'e inequality holds with the optimal constant $\frac{N-1}{KN}$ instead of the constant $\frac{1}{K}$.

\subsubsection{Transportation inequality}

Another family of functional inequalities which is strongly related to Ricci curvature bounds is the family of transportation inequalities. Here, we focus on transportation inequalities $T_2$ , which control the Wasserstein distance $W_2(\mu,\nu)$ between two probability measures by the square root of the relative entropy $\sqrt{H_\nu(\mu)}$. 

\begin{defi}
A measured length space is said to satisfy the transportation inequality $T_2(K)$, where $K>0$ if, for every $\mu$ in $P_2(X)$, \begin{equation}W_2(\mu,\nu)^2 \leq \frac{2}{K} H_\nu(\mu).\end{equation}
\end{defi}

Transportation inequalities have been first studied by Talagrand in \cite{talagrand1996transportation}, which established transportation inequality $T_2$ for gaussian measures, and proved tensorization properties. In \cite{OttoVill}, Otto and Villani derive the inequality $T_2(K)$ from $LSI(K)$ in the weighted Riemannian setting, and conversely prove that the inequality $T_2(K)$ yield a logarithmic Sobolev inequality, with some degradation on the constant. The implication $LSI(K) \Rightarrow T_2(K)$ has been extended to the length space setting by Lott and Villani in \cite[Theorem 6.22]{LottVill}, under some regularity assumption on the ambient space. In the same paper, the authors prove directly the implication $CD(K,\infty) \Rightarrow T_2(K)$.

\medskip

 In this subsection, we prove that the condition $PL(K)$, with $K>0$, yield a transportation inequality.
\begin{thm}\label{Talagrand}
Let $(X,d,\nu)$  be a measured length space satisfying $PL(K)$ for some $K>0$. Then the transportation inequality $T_2(K)$ holds true.
\end{thm}

This fact is already known in smoother cases. For example, it is stated and proven in a weighted Euclidean case in \cite{BobkovLedoux}.
The proof of Theorem \ref{Talagrand} in a length space setting follows the same lines as in the Euclidean one, and uses the following proposition, which is the variational formulation of the inequality.

\begin{prop}\label{TalagrandVariational}
The classical formulation of the Talagrand inequality: \begin{equation}\label{TalagrandEq}\forall \mu \in P(X), \ W_2(\mu,\nu)^2 \leq \frac{2}{K} H_\nu(\mu)\end{equation} is equivalent to its variational formulation: \begin{equation}\label{TalagrandVarietionalEq}\forall g \in \mathcal{C}(X), \ \int_X e^{K Q_1g(x)}d\nu \leq e^{K \int_X g d\nu(x)},\end{equation} where: $$Q_1g(x) := \inf_{y \in X} \left( g(y) +\frac{d(x,y)^2}{2} \right).$$
\end{prop}

The notation $Q_1g$ comes from the Hamilton-Jacobi semi-group. See \cite{BobvovGentilLedoux} for details and a proof of this proposition.

\medskip

\textbf{Proof of Theorem \ref{Talagrand}:} For $g$ be in $\mathcal{C}(X)$, and $t \in [0,1]$, define $$u_t(x) := \exp\left(-t K g(x) \right),$$ and $$v_t(y) := \exp\left( (1-t) K Q_1g(y) \right).$$ Then
\begin{eqnarray*}
u_t(x)^{1-t}v_t(y)^{t} &=& \exp\left(-Kt(1-t) (g(x)-Q_1g(y))\right) \\
&\leq& \exp\left( -Kt(1-t) \left(g(x)-\left(g(x)-\frac{d(x,y)^2}{2}\right)\right) \right) \\
&=& \exp\left(Kt(1-t) \frac{d(x,y)^2}{2} \right).
\end{eqnarray*}
Therefore $$\forall x,y \in X, \ u_t(x)^{1-t}v_t(y)^t \exp\left(-Kt(1-t) \frac{d(x,y)^2}{2}\right) \leq 1.$$
Applying the Pr\'ekopa-Leindler inequality to the triple of functions $(u_t,v_t,1)$ yields $$\left(\int_X u_t d\nu \right)^{1-t}\left( \int_X v_t \right)^t \leq 1.$$
Thus $$\int_X \exp\left( (1-t) K Q_1g(y) \right) d\nu(y) \leq \left(\int_X \exp\left(-t K g(x)\right)d\nu(x)\right)^{-\frac{1-t}{t}}.$$
Since $$\lim_{t \rightarrow 0} \int_X \exp\left( (1-t) K Q_1g(y) \right) d\nu(y) = \int_X \exp\left( K Q_1g(y) \right) d\nu(y), $$ it suffices to show that $$\limsup_{t \rightarrow 0} \left(\int_X \exp\left(-t K g(x)\right)d\nu(x)\right)^{-\frac{1-t}{t}} \leq e^{K \int_X g d\nu(x)}.$$
And this indeed holds as a consequence of Jensen's inequality
\begin{eqnarray*}
log\left(\int_X \exp\left(-t K g(x)\right)d\nu(x)\right)^{-\frac{1-t}{t}} &=& -\frac{1-t}{t} log \int_X \exp\left(-t K g(x)\right)d\nu(x) \\
&\leq& -\frac{1-t}{t} \int_X -tKg(x) d\nu(x) \\
&=& (1-t) K \int_X g(x) d\nu (x). 
\end {eqnarray*} This proves the theorem.$\square$

\subsection{Links with the Brunn-Minkowski inequality}

The Pr\'ekopa-Leindler inequality is known to be the functional version of the multiplicative Brunn-Minkowski inequality.

\begin{defi}
A measured length space $(X,d,\nu)$ is said to satisfy the mutliplicative Brunn-Minkowski inequality, or $BM$, if, for every measurable sets $A$ and $B$, and every $t$ in $[0,1]$: $$\nu^*(Z_t(A,B)) \geq \nu(A)^{1-t} \nu(B)^t.$$
\end{defi}

The notation $\nu^*$ corresponds to the outer measure associated to $\nu$, and is useful to avoid measurability issues. However, if $(X,d)$ is compact, then the intermediate sets $Z_t(A,B)$ are automatically measurable.\\
It is very easy to deduce $BM$ from $PL(0)$, by taking indicator functions of the sets $A$ and $B$. This subsection is about the converse of this fact:

\begin{thm}\label{BMPL}
Let $(X,d,\nu)$ be a length space satisfying $BM$. Then it also satisfies $PL(0)$.
\end{thm}

The proof of this result requires the following Lemma, known as the ``multiplicative Pr\'ekopa-Leindler inequality''.

\begin{lem}
Let $f,g,h$ in $L^1(\mathbb{R}_+,\mathbb{R}_+)$ satisfying $$\forall u,v \in \mathbb{R}_+, \ h(u^{1-t}v^t) \geq f(u)^{1-t}g(v)^t.$$ Then: $$\int_{\mathbb{R}_+} h(w) dw \geq \left(\int_{\mathbb{R}_+} f(u) du\right)^{1-t} \left(\int_{\mathbb{R}_+} g(v) dv\right)^t$$
\end{lem} A proof of this Lemma can be found in \cite{Uhrin94}.

\bigskip

\textbf{Proof of Theorem \ref{BMPL}:} Let $f,g,h : X \rightarrow \mathbb{R}_+$ be a triple of functions satisfying the hypothesis of the Pr\'ekopa-Leindler inequality. The idea consists in studying the level sets of these functions. They satisfy the inclusion $$\forall t \in [0,1], \ \forall u,v \in \mathbb{R}_+, \ Z_t(\{f \geq u\},\{g \geq v\}) \subset \{h \geq u^{1-t}v^t\}.$$ So the $BM$ inequality yields $$\nu\left(\{f \geq u\}\right)^{1-t} \nu(\{g \geq v\})^t \leq \nu(\{h \geq u^{1-t}v^t\}).$$
Define the function $\overline{f} : \mathbb{R}_+ \rightarrow \mathbb{R}$ by: $\overline{f}(u) = \nu(\{f\geq u\})$, and $\overline{g}$, $\overline{h}$ in the same way. Fubini's theorem indicates that  $\int_{\mathbb{R}_+} \overline{f}(u) du = \int_X f d\nu$, and so on for $\overline{g}$ and $\overline{h}$.
It is clear that the triple $\overline{f} ,\overline{g} ,\overline{h}$ satisfies the hypothesis of mulitplicative Pr\'ekopa-Leindler inequality. Consequently, $$\int_{\mathbb{R}_+} \overline{h}(w) dw \geq \left(\int_{\mathbb{R}_+} \overline{f}(u) du\right)^{1-t} \left(\int_{\mathbb{R}_+} \overline{g}(v) dv\right)^t.$$ This proves Theorem \ref{BMPL}. $\square$

\medskip

Theorem \ref{BMPL}, combined with Theorem \ref{PLRic} and \cite[Theorem 7.3]{LottVill}, has a nice interpretation. On a compact Riemannian manifold, the three conditions $BM$, $PL(0)$, $CD(0,\infty)$ are all equivalent to the classical Ricci curvature bound condition $\ric\geq0$.

\section{The Borell-Brascamp-Lieb inequality $BBL(0,N)$}

In this section, we investigate some of the properties of the $BBL(0,N)$ condition. The main goal is to show that the $BBL(0,N)$ condition is a natural generalization of the couple of conditions $\text{Ric} \geq 0$ and $\textrm{dim} \leq N$ on a Riemannian manifold, by following the same lines as in the last section. We first focus on the coherence with the Riemannian case. Then we prove a tensorization property, similar to Theorem \ref{PLTenso}. We finally give some consequences of the $BBL(0,N)$ condition: the Bishop-Gromov theorem, that allows us to control the growth of volume of balls around a given point (Theorem \ref{BishopGromov}), and yields a bound on the Hausdorff dimension of the underlying metric space (Corollary \ref{HausDim}). With an additional $PL(K)$ condition with $K>0$, and a compacity assumption, the Bishop-Gromov theorem also yields a control of the diameter of the space (Theorem \ref{BishopGromov}).

\subsection{The Riemannian case}

In \cite{CMS01}, Cordero-Erausquin, McCann and Schmuckenschl\"ager proved that the $BBL(0,N)$ inequality holds on any Riemannian manifold $(M,g,d\textrm{vol})$ satisfying $\ric \geq 0$ and $\textrm{dim}(M) \leq N$. We establish here a converse statement in the weighted Riemannian case, which can be seen as the ``$BBL(0,N)$ analogous'' to Theorem \ref{PLRic}.

\begin{thm}\label{BBLRic}
Let $(M,g, e^{-V}d\textrm{vol})$ be a Riemannian manifold satisfying the condition $BBL(0,N)$ for a certain $N>0$. Then $\ric_m+\hess_mV \geq 0$ at each point $m$ of the manifold and $\textrm{dim}(M) \leq N$.
\end{thm}

\textbf{Proof :} Recall that the $BBL(0,N)$ condition implies $PL(0)$. Then, by Theorem \ref{PLRic}, we have immediatly $\ric_m+\hess_mV \geq 0$. The inequality $\textrm{dim}(M) \leq N$ is a consequence of Corollary \ref{HausDim}, because the dimension of a smooth manifold coincide with its Hausdorff dimension. $\square$

\subsection{Tensorization}

In this paragraph, we establish a tensorization property for the $BBL(0,N)$ inequality, which is very similar to the tensorization property found for the $PL(K)$ condition (Theorem \ref{PLTenso}). However, in this setting there is no particular choice of the metric on the product space.

\begin{thm}\label{BBLTenso}
Let $(X_1,d_1,\nu_1)$ and $(X_2,d_2,\nu_2)$ be two measured length spaces satisfying respectively $BBL(0,N_1)$ and $BBL(0,N_2)$, for $N_1,N_2 \geq 1$. 
Let $d$ be a distance on $X := (X_1 \times X_2)$, $\nu := \nu_1 \otimes \nu_2$ and $N := N_1+N_2$.
Then $(X,d,\nu)$ satisfies $BBL(0,N)$.
\end{thm}

The fact that $N = N_1+N_2$ is coherent with the intuitive meaning of $N$ as an upper bound on the dimension of the space. This intuition will be made rigorous at Corollary \ref{HausDim}.

\medskip

\textbf{Proof:} Let $p\geq -\frac{1}{N}$. In particular, $p \geq -\frac{1}{N_1}$. Let $t\in [0,1]$ and $f,g,h : X \rightarrow \mathbb{R}_+$ satisfying $$\forall x,y \in X, \ \forall z \in Z_t(x,y), \ h(z) \geq \mathcal{M}_t^p(f(x),g(y)).$$
Fix $x_2,y_2,z_2 \in X_2$ such that $z_2 \in Z_t(x_2,y_2)$. Then \begin{equation}\forall x_1,y_1 \in X_1, \ \forall z_1 \in Z_t(x_1,y_1), \ (z_1,z_2) \in Z_t((x_1,y_1),(x_2,y_2)).\end{equation} In particular, the functions $f_{x_2},g_{y_2},h_{z_2} : X_1 \rightarrow \mathbb{R}_+$ satisfy $$\forall x_1,y_1 \in X_1, \ \forall z_1 \in Z_t(x_1,y_1), \ h_{z_2}(z_1) \geq \mathcal{M}_t^p(f_{x_2}(x_1),g_{y_2}(y_1)).$$ Applying the $BBL(0,N_1)$ inequality on $(X_1,d_1,\nu_1)$ yields \begin{equation}\int_{X_1} h_{z_2} d\nu_2 \geq \mathcal{M}_t^{p/(1+N_1 p)}\left(\int_{X_1}f_{x_2} d\nu_1, \int_{X_1} g_{y_2} d\nu_1\right).\end{equation} Set $p'=p/(1+N_1 p)$, and define the function $\tilde{f} : X_2 \rightarrow \mathbb{R}_+$ by $\tilde{f}(x_2) := \int_{X_1} f_{x_2} d\nu_1$ and the functions $\tilde{g}$ and $\tilde {h}$ in the same way. It is clear that this triple of functions satisfy the hypotheses of the $BBL_{p'}(0,N_2)$ inequality on $(X_2,d_2,\nu_2)$. Consequently, \begin{equation}\int_{X_2} \tilde{f} d\nu_2 \geq \mathcal{M}_t^{p'/(1+N_2 p')} \left( \int_{X_2} \tilde{f} d\nu_2 , \int_{X_2} \tilde{g} d\nu_2 \right).\end{equation} Immediate properties of product measures tell that $$\int_{X_2} \tilde{f} d\nu_2 = \int_{X_1 \times X_2} f d(\nu_1 \otimes \nu_2),$$ and similarly for $\tilde{g}$ and $\tilde{h}$. Moreover, it is easy to see that $$\frac{p'}{1+N_2 p'} = \frac{p}{1+N_1 p} \cdot \frac{1}{1+N_2 \frac{p}{1+N_1 p}} = \frac{p}{1+(N_1+N_2)p}.$$ Putting everything together gives $$\int_X h d\nu \geq \mathcal{M}_p^t \left(\int_X f d\nu, \int_X g d\nu \right),$$ and the theorem is proven. $\square$

\medskip

As the $BBL(0,N)$ condition implies the $PL(0)$ condition, it is possible to improve the statement of Theorem \ref{PLTenso} in the case where $K=0$: this Theorem holds true for every metric on the product space, and not only for the Euclidean product metric.

\subsection{Geometric consequences}

The Bishop-Gromov inequality is one of the most famous results in Riemannian comparison geometry. It gives a control of the growth of the volume of balls around a fixed point on a Riemannian manifold with bounded Ricci curvature. The reader can find a proof of this theorem in any reference book on Riemannian geometry, for example \cite{GHL}.

\begin{thm}
Let $(M,g,d\textrm{vol})$ be a $n$-dimensional Riemannian manifold satisfying $\ric \geq K$.
Let us denote by $V(r)$ the volume of the ball of radius $r$ in $S_K^n$, where $S^n_K$ is a complete simply connected $n$-dimensional Riemannian manifold of constant sectional curvature $K$. Then the quotient $$\frac{\textrm{vol}(B_p(r))}{V(r)}$$ is non-increasing in $r$.
\end{thm}

In the case where $K=0$, $S^n_K$ is just the Euclidean space $\mathbb{R}^n$, and the function $V(r)$ is proportional to $r^n$.
The goal of this paragraph is to show that the same control is still valid for any measured length space satisfying the $BBL(0,N)$ condition:

\begin{thm}\label{BishopGromov}
Let $(X,d,\nu)$ be a measured length space satisfying the $BBL(0,N)$ condition for a certain $N>0$. Let $x_0 \in X$, and $0<r_1\leq r_2$. Then: \begin{equation}\label{BishopGromovEq}\nu(B_{r_2}(x_0)) \leq \left(\frac{r_2}{r_1}\right)^N \nu(B_{r_1}(x_0)).\end{equation}
\end{thm}

\medskip

\textbf{Proof:} Set $\eta > 0$, $\varepsilon > 0$ small enough to have $\varepsilon < \eta$ and $t:=(r_1+\varepsilon)/r_2 \in [0,1]$, and $$u(x) = 1_{B_{\varepsilon}(x_0)}(x), \ v(y) = 1_{B_{r_2}(x_0)}(y), \ w(z) = 1_{B_{r_1+\eta}(x_0)}(z).$$ 

The first step is to prove that the triple of functions $u,v,w$ satisfy the hypothesis of the $BBL(0,N)$ inequality \begin{equation}\label{BBLHypo}\forall x,y \in X, \forall z \in Z_t(x,y), \ \forall p \geq -\frac{1}{N}, \ w(z) \geq \mathcal{M}_t^p(u(x),v(y)).\end{equation} As $\mathcal{M}_t^p(u(x),v(y)) \leq 1$, it suffices to prove \ref{BBLHypo} for any $z \in X$ such that $w(z) = 0$, i.e. for $z \notin B_{r_1+\eta}(x_0)$.
Let $z \in X$ such that $z \notin B_{r_1+\eta}(x_0)$, and $x,y \in X$ such that $z \in Z_t(x,y)$. The formula (\ref{BBLHypo}) will be proven if either $u(x)=0$ or $v(y)=0$. Suppose that $u(x) \neq 0$ and $v(y) \neq 0$. This means that $d(x,x_0) \leq \varepsilon$ and $d(y,x_0) \leq r_2$. Then 
\begin{eqnarray*}
d(x_0,z) &\leq& d(x_0,x) + d(x,z) \\
&\leq& \varepsilon + t d(x,y) \\
&\leq&  \varepsilon + (\frac{r_1+\varepsilon}{r_2})(d(x,x_0)+d(x_0,y)) \\
&\leq& \varepsilon + (r_1+\varepsilon) (1+\frac{\varepsilon}{r_2}) \\
&\leq& r_1+\eta.
\end{eqnarray*}
This contradicts the fact that $z \notin B_{r_1+\eta}(x_0)$, so the formula (\ref{BBLHypo}) is proven.\\
Now, the $BBL(0,N)$ inequality implies $$\forall p  \geq -\frac{1}{N}, \ \int_X w d\nu \geq \mathcal{M}_t^{p/(1+Np)}\left(\int_X u d\nu, \int_X vd\nu \right).$$ This could be written as $$\forall p  \geq -\frac{1}{N}, \ \nu(B_{r_1 + \eta}(x_0)) \leq t^{\frac{1+Np}{p}} \nu(B_{r_2}(x_0))$$ 
As $t=\frac{r_1}{r_2}$, and $\frac{1+Np}{p}=\frac{1}{p}+N$, $$\nu(B_{r_1+\eta}(x_0)) \leq \left(\frac{r_1+\varepsilon}{r_2}\right)^{N+1/p}\nu(B_{r_2}(x_0)).$$ The theorem is proven by letting first $p$ go to $+\infty$ and then by letting $\eta$ (and so $\varepsilon$) go to $0$. $\square$

\medskip

A first corollary of Theorem \ref{BishopGromov} allows for a control of the Hausdorff dimension of a measured length space satisfying the $BBL(0,N)$ condition, which gives a rigorous explanation of the meaning of the parameter $N$ as an upper bound on the dimension of $(X,d)$.A proof of this fact can be found in \cite{Sturm06II}.

\begin{cor}\label{HausDim}
Let $(X,d,\nu)$ be a measured length space satisfying the $BBL(0,N)$ condition. Then the support of $\nu$ has Hausdorff dimension $\leq N$.
\end{cor}

\medskip

We finish this paragraph by giving a corollary of both Theorems \ref{Talagrand} and \ref{BishopGromov}, which allows to control the diameter estimate of a compact length space satisfying both conditions $PL(K)$ and $BBL(0,N)$, where $K>0$ and $N\geq 1$. This kind of estimate is known in the Riemannian setting as Bonnet-Myers Theorem.

\begin{thm}\label{BonnetMyersRiem}
Let $(M,g)$ be a complete $N$-dimensional Riemannian manifold satisfying $\ric \geq K$, where $K>0$. Then $(M,g)$ is compact and $$\diam(M,g) \leq \pi \sqrt{\frac{N-1}{K}} = \diam\left(\mathbb{S}^N_{\frac{1}{\sqrt{K}}}\right).$$
\end{thm}

An analogous theorem in a length space setting is stated as follows.

\begin{thm}\label{BonnetMyers}
Let $(X,d,\nu)$ be a compact measured length space satisfying both conditions $PL(K)$, for a certain $K>0$, and $BBL(0,N)$, for a certain $N>0$. Then the diameter of $X$ satisfies: \begin{equation}\label{BonnetMyersEq}\diam(X) \leq C \sqrt{\frac{N}{K}},\end{equation} where $C$ is a universal constant (C=7.7 is a good choice).
\end{thm}

The length space version of Bonnet-Myers is weaker than its Riemannian version because Theorem \ref{BonnetMyersRiem} does not assume a priori the compacity of $(M,g)$, and of course because the right-hand side of the inequality is more precise.

\medskip

\textbf{Proof of Theorem \ref{BonnetMyers}:} In \cite{LottVill}, it was proven that the inequality \eqref{BonnetMyersEq} is satisfied whenever both Talagrand inequality $T(K)$ (inequality \eqref{TalagrandEq}) and Bishop-Gromov inequality (inequality ~\eqref{BishopGromovEq}) with constant $N$ are true. But inequality \eqref{TalagrandEq} is a consequence of Theorem \ref{Talagrand}, and inequality \eqref{BishopGromovEq} is a consequence of Theorem \ref{BishopGromov}. $\square$

\bibliography{prekobib}
\bibliographystyle{alpha}

\end{document}